\documentclass[11pt, reqno]{glm-t}
\usepackage{amsthm, amscd, amsfonts,  graphicx, color, xcolor}
\usepackage{amssymb,mathrsfs,amsthm,amsmath,makeidx,verbatim, euscript,eufrak}
\usepackage{geometry}
\usepackage{amssymb,mathrsfs,amsthm}
\usepackage{multirow}
\usepackage[english]{babel}
\usepackage[all]{xy}
\usepackage{mathtools}
\newtheorem{thm}{Theorem}[section]
\newtheorem{prop}[thm]{Proposition}
\newtheorem{lem}[thm]{Lemma}
\newtheorem{coro}[thm]{Corollary}
\theoremstyle{definition}
\newtheorem{defx}[thm]{Definition}
\newtheorem{rem}[thm]{Remark}

\newcommand{\R}{I\!\!R}

\begin{document}

\title{{\color{green!42!blue}\LARGE{\bf Partial Dynamical Systems of $L^p$-Spaces and their Stability Spaces}}}

\author[a1]{N. O. Okeke}

\author[a2]{M. E. Egwe$^*$\corref{*}}


\address[a1]{ Physical and Mathematical Sciences, Dominican University, Ibadan.\\$nickokeke@dui.edu.ng$  \vskip 0.1cm}

\address[a2]{Department of Mathematics, University of Ibadan, Ibadan, Nigeria.\\$murphy.egwe@ui.edu.ng$ {\tt}}

\cortext[c1]{M. E. Egwe}
\setcounter{page}{1}
\vol{15 (2025),\; Issue 2} \pages{65-73}


\recivedat{received date}

\authors{N.O. Okeke, M.E.Egwe}

\doi{\href{doi.org/10.31559/glm}{}}

\begin{abstract}
In this paper, we consider the partial dynamical systems of the locally convex, $L^p(\Omega)-$ spaces defined by the action of the smooth algebra $\mathscr{K}(\Omega)$ through its nets. Slice analysis is then employed to show that the Sobolev spaces $W^{k,p}(\Omega)$ are the stable states or space of these partial dynamical systems as limit spaces of the convolution actions of the smooth algebra $K(\Omega)$ on the Banach spaces $L^p(\Omega)$. Thus, the Sobolev spaces $W^{k,p}(\Omega)$ are closed subspaces of the $Lp(\Omega)$-spaces under convolution product and weak derivatives, with the weak derivative operators acting as equivariant maps of the slice spaces. \copyright 2025 All rights reserved.
\end{abstract}

\maketitle



\section{Introduction}
Given the differential equation of the form $f'(x) = kf(x)$, where $k$ is a real constant,  solutions of the form $f(x) = f(0)e^{kx}$ on a neighbourhood $U_{x_o}$ of a fixed $x_o$ such that $f(x) > 0 \; \forall \; x \in U_{x_o}$ can be established. In its simplest form, the differential equation shows an exponential transform of a connected component or a line segment. The differential equation has solution everywhere since all continuous functions $f: [a,b] \to \R$ satisfying the condition $|f'(x)| \leq k|f(x)|$ are entire functions on $\R$. If the domain is $\Omega \subset \R^n$, the idea of parametrization is used to obtain a differential equation of the form $f'(t) = X(f(t))$ with initial condition $f(0) = x_o$. What determines the solution type is the nature of the differential operator $X$ acting on the space of parameterized functions $x = f(t)$.

To move from the simple differential equation given on $\R$ to the parameterized form which is its embedding in the manifold $\R^n$, where $X$ is a Lipschitz vector field giving a unique solution for every $x_o$ in $U \subset \R^n$, the idea of diffeomorphism is used. Diffeomorphisms between open subsets of $U$ are generated by setting $\phi_t(x_o) = f(t)$. The intention is to study the dynamics of a system of these differential equations using the weak solutions which are globalizable. Since these differential equations are flows in the manifold, their generalization would require smooth vector fields as generators.

Because the condition of differentiability adds uniformity to continuity of functions, one can generalize convergence of sequences to nets. Hence, the concept of a Cauchy sequence is generalized to that of convergent net. This locates solution space of this form of differential equations within the Banach $L^p$-spaces. To analyse the solution space, we use the idea of smooth kernels or unitary multipliers, which are generalized from the definition of a continuous net $\{\rho_\varepsilon\}_{\varepsilon > 0}$ in \cite{Jost2005}, as a directed set of continuous nonnegative functions $\rho_\varepsilon$ on $\R^n$ with supports $spt(\rho_\varepsilon)$ in $B(\varepsilon)$ and $\displaystyle \int \rho_\varepsilon d\mu^n = 1$, where $\mu^n$ is the Lebesgue measure on $\R^n$.

The concept of a smooth net $\phi_\varepsilon$ follows from continuous net when each component is non-negative, uniformly bounded and continuous, integrable, and infinitely differentiable. That is, each $\phi$ is of the class $C^\infty(\R^n)$ supported on a compact subset of $\R^n$. A smooth algebra is defined on a set of smooth nets as follows.
\begin{defx}
A smooth algebra is a collection of smoothing nets defined on the open subset $\Omega$ having the following properties (i) $\phi_\varepsilon \in C^\infty_o(\R^n)$, (ii) $\text{supp}(\phi_\varepsilon) \subset \bar{B}(0,\varepsilon)$, and (iii) $\displaystyle \int \phi_\varepsilon(x)dx = 1$. It is a smooth algebra of units  denoted by $\mathscr{K}(\Omega)$.
\end{defx}
Any smooth function $\phi : \R^n \to [0,\infty)$ supported on a unit ball $B(1)$, which can be supported on any closed ball by proper modification such that $\displaystyle \int \phi d\mu^n = 1$, defines an approximate identity by \[\phi_\varepsilon(x) = \varepsilon^{-n}\phi(\frac{x}{\varepsilon}). \eqno{(1.0)}  \]
Examples of such smooth kernels are given in \cite{AlGwaiz92}, \cite{Jost2005} and \cite{Grosseretal2001}.
\begin{rem}
The smooth net is different from the delta net and the product nets as defined and explained in \cite{Grosseretal2001}. It is a natural tool intrinsically related to the local dynamical system on the Banach subspace $L^p(\Omega)$, which is defined by its action. The action results in fibre structure which carries a distributional system and a weak differential structure. The employment of the smooth net in the canonical definition of the generalized algebra of distributions in section 3.3.2 of \cite{Grosseretal2001} points to its importance as intrinsic tool for analysis on smooth manifolds.
\end{rem}
\section{Distributional Dynamical Systems}
 the smooth algebra acts by its nets on the Banach spaces $L^p(\Omega)$. This convolution action will be used to study the distributional solution of a differential equation constituting a distributional (partial) dynamical system caused by the action. The smooth algebra $\mathscr{K}(\Omega)$ defines a convolution action on $L^p(\Omega)$ as follows.
\begin{defx}
Given the smooth algebra $\mathscr{K}(\Omega)$, we defined its convolution action on $L^p(\Omega)$ via its net $(\phi_\varepsilon)$ as follows \begin{equation} \mathscr{K}(\Omega) \times L^p(\Omega) \to L^p(\Omega), \phi_\varepsilon(f) := \frac{1}{\varepsilon^n}\int_\Omega \phi(\frac{x-y}{\varepsilon})f(y)dy = f_\varepsilon(x).\end{equation}
\end{defx}
The denseness of $C^\infty(\Omega)$ in $L^p(\Omega)$ means the smoothing algebra action is closed on a subspace of $L^p(\Omega)$.  A fibre structure follows from this action because it defines a local distributional system on the subspace $L^p(\Omega')$ for every compact subset $\Omega' \subset \subset \Omega$. Subsequently, given that $C^\infty_o(\Omega) = \{\phi \in C^\infty(\Omega) : supp(\phi) \subset \subset \Omega\}$, the following result shows the denseness of $C^\infty_o(\Omega)$.
\begin{prop}\normalfont
Given $\Omega' \subset \subset \Omega$ and $\varepsilon < dist(x, \partial \Omega)$ then $\phi_\varepsilon(f) = f_\varepsilon$ is an orbit of the $\mathscr{K}(\Omega)$-action on $L^p(\Omega)$.
\end{prop}
\begin{proof}
From the definition of the convolution action above, it shows $f_\varepsilon$ to be in $C^\infty(\Omega')$ once $\phi_\varepsilon$ is in $C^\infty(\Omega')$, for every compact subset $\Omega'$ in $\Omega$. Taking $D^k$ to be all the $k$-th partial derivatives, by the uniform boundedness and continuity of the elements of the net, there exists a constant $c_k$ such that $\underset{x \in \R^n}\sup|D^k\phi_\varepsilon(x)| \leq c_k$. This is considered in \cite{Grosseretal2001} as $c$-boundedness of all $\phi_\varepsilon$ in any compact subset. Thus, $(f_\varepsilon) \subset C^\infty(\Omega')$ is a net of smooth flows converging to $f \in L^p(\Omega)$; and is a representation of a connected component (or a net) of $\mathscr{K}(\Omega)$ on the space $\text{Diff}(L^p(\Omega))$ of diffeomorphisms of $L^p(\Omega)$.
\end{proof}
Subsequently, since the action of the smooth algebra $\mathscr{K}(\Omega)$ by its nets is by convolution product operation on the $L^p$-spaces or their subspaces and commutes with weak derivative, the convolution product and weak derivative are used to study the Sobolev spaces. The weak derivative introduces the Sobolev spaces as its invariant (or closed) subspaces with a weak differential structure regularized by the smooth action. They are subspaces of the $L^p$-spaces (Lebesgue integrable functions) possessing \emph{weak derivatives}.

Because continuous differentiability is a necessary condition for partial differentiation with respect to all variables, it is a requirement for definition of weak derivatives and the generalization of the definition to higher dimensional spaces as given in section 20-4 of \cite{Jost2005} using the idea of partial differentiation.
\begin{defx}
Denote the weak derivative of $f$ in $x^i$-direction with $D_if$, and $Df = (D_1f, \cdots, D_nf)$ when $f$ has a weak derivative for all $i = 1,\cdots, n$. Then with the multi-index $\alpha := (\alpha_1, \cdots, \alpha_n), \alpha_i \geq 0, (i = 1, \cdots, n)$, and the norm given as $ |\alpha| := \overset{n}{\underset{i=1}\sum}\alpha_i > 0$, we have \[ D_\alpha \varphi := \left(\frac{\partial}{\partial x^1}\right)^{\alpha_1} \cdots \left(\frac{\partial}{\partial x^n}\right)^{\alpha_n}\varphi = \frac{\partial^{|\alpha|}\varphi}{\partial x_1^{\alpha_1} \cdots \partial x_n^{\alpha_n}} \; \text{for } \varphi \in C^{|\alpha|}(\Omega). \]
Then a function $u \in L^1_{loc}(\Omega)$ is called the $\alpha$-th weak derivative of $f$ if \[ \int_\Omega \varphi udx = (-1)^{|\alpha|}\int_\Omega f D_\alpha\varphi dx, \; \text{for all } \varphi \in C^\infty_c(\Omega) \subset C^{|\alpha|}_o(\Omega). \]
\end{defx}
\begin{rem}
The following remarks clarify the definition.\\
(1) Weak derivatives are unique up to a set of measure zero. Hence, the weak derivative of $f$ is also written $D_\alpha f$. We therefore have $\langle D_\alpha f, \varphi \rangle = (-1)^{|\alpha|}\langle f, \nabla \varphi \rangle$.\\
(2) Since all functions with weak derivatives are not all in $C^1(\Omega)$, and not all functions in $L^1_{loc}(\Omega) \subset L^p(\Omega)$ have weak derivatives, special spaces for all the functions with weak derivatives are defined as Sobolev spaces as follows in \cite{Jost2005}, 20.5.
\end{rem}
\begin{defx}
For $k \in \mathbb{N}, 1 \leq p \leq \infty$, the Sobolev spaces are defined as the sets $W^{k,p}(\Omega) := \{f \in L^p(\Omega) : D_\alpha f \; \text{exists and is in } L^p(\Omega) \; \forall \; |\alpha| \leq k \}$. The norm on the spaces is defined as \[||f||_{W^{k,p}(\Omega)} = ||f||_{k,p} := \left(\underset{|\alpha| \leq k}\sum \int_\Omega |D_\alpha f|^p \right)^{\frac{1}{p}} \; \text{for } 1 \leq p < \infty. \] and \[||f||_{W^{k,\infty}(\Omega)} = ||f||_{k,\infty} := \underset{|\alpha| \leq k}\sum \underset{x \in \Omega}{\text{ess}\sup}|D_\alpha f(x)|. \]
\end{defx}
The definition of weak derivative is compatible with the convolution product. The Sobolev spaces are subspaces of $L^p(\Omega)$ spaces closed under convolution product and weak derivative. It follows that the action of the smoothing algebra $\mathscr{K}(\Omega)$ on $L^p(\Omega)$ commutes with (weak) differential operator on Sobolev spaces. This is given in \cite{Jost2005}, 20.6 as follows.
\begin{prop}\normalfont
Let $f \in L^1_{loc}(\Omega)$ and assume that $D_if$ exists. If $\varepsilon < dist(x,\partial \Omega)$ then the convolution action of the smooth algebra $\mathscr{K}(\Omega)$ by its smooth net $\phi_\varepsilon(f) = f_\varepsilon$ commutes with the (weak) differential operator $D_i$.
\end{prop}
\begin{proof}
The commutativity of the $\mathscr{K}(\Omega)$-action with the weak/partial differential operator $D_i$ given as $(D_i \star \phi_\varepsilon)(f) = D_i(\phi_\varepsilon(f)) = D_if_\varepsilon = \phi_\varepsilon(D_if) = (\phi_\varepsilon \star D_i)(f)$ implies that the weak (partial) differential operators $D_i$ are $\mathscr{K}(\Omega)$ invariant operators on the Sobolev spaces $W^{k,p}(\Omega)$ which are dense subspaces of $L^p(\Omega)$. Thus, the Sobolev spaces $W^{k,p}(\Omega)$ are $\mathscr{K}(\Omega)$-invariant spaces. This means that the smooth algebra $\mathscr{K}(\Omega)$ encode the local symmetry of the Sobolev spaces $W^{k,p}(\Omega)$.
\end{proof}
This is expressed in the following result (cf. 20.7 of \cite{Jost2005}).
\begin{thm}\normalfont
The weak differential operator $D_\alpha : L^p(\Omega) \to L^p(\Omega)$ is closed and $\mathscr{K}(\Omega)$-equivariant transformation on the Sobolev space $W^{k,p}(\Omega)$.
\end{thm}
\begin{proof}
Let $\phi_\varepsilon(f) = f_\varepsilon$ be the action of $\mathscr{K}(\Omega)$, then $\phi_\varepsilon(u) = u_\varepsilon = D_i(f_\varepsilon) = D_i(\phi_\varepsilon(f))$. Then as $\varepsilon \to 0$ $f_\varepsilon \to f$ in $L^p(\Omega)$, and $u_\varepsilon \to u$ in $L^p(\Omega')$ for any $\Omega' \subset \subset \Omega$. On the other hand, suppose that $\phi \in C^1_o(\Omega)$. Then $D_i\phi, (i = 1, \cdots, n)$ are bounded; and the dominated convergence theorem guarantees that \[\int_\Omega f(x)\frac{\partial}{\partial x^i}\phi(x)dx = \underset{\varepsilon \to 0}\lim \int_\Omega f_\varepsilon(x)\frac{\partial}{\partial x^i}\phi(x)dx \] \[ = -\underset{\varepsilon \to 0}\lim \int D_if_\varepsilon(x) \phi(x)dx = - \int_\Omega u(x)\phi(x)dx. \] Thus, $u = D_if$, and $f, f_\varepsilon, D_i(f_\varepsilon), D_if \in W^{k,p}(\Omega)$ (by definition).
\end{proof}
\begin{lem}\normalfont
The smoothing action of $\mathscr{K}(\Omega)$ defines a partial dynamical system on $L^1_{loc}(\Omega)$ by its nets.
\end{lem}
\begin{proof}
Given the definition of partial dynamical system in \cite{Exel2017}, (which has to do with the extension of the solution to the initial value problem $f'(t) = X(f(t)), \; f(0) = x_o$ on some open interval about zero in $\R^n$), a partial dynamical system is defined by composition of diffeomorphisms $\phi_t(x_o) = f(t)$ between open subsets $U \subset \R^n$, where $x_o \in U$. Thus, the implication \begin{equation} \phi_{s+t}(x) = \phi_s(\phi_t(x)) \; \implies \; \phi_s \circ \phi_t \subseteq \phi_{s+t}, \end{equation}
means $\phi_{s+t}$ extends $\phi_s \circ \phi_t$.

Applying this to the space of functions, we have a partial action of the smooth algebra $\mathscr{K}(\Omega)$ as a pair $(\{W_\varepsilon\},\{\phi_\varepsilon\})$ where $\phi_\varepsilon$ are diffeomorphisms defined on the subsets $W_\varepsilon$ of $L^p(\Omega)$. A partial dynamical system is therefore given by the quadruple $(L^p(\Omega), \mathscr{K}(\Omega), W^{k,p}(\Omega), \phi_\varepsilon)$. Thus the net $\phi_\varepsilon \subset \mathscr{K}(\Omega)$ acting on $W^{k,p}(\Omega) \subset L^p(\Omega)$ is a partial dynamical system on the $L^p$-spaces.
\end{proof}
Though the smooth algebra $\mathscr{K}(\Omega)$ defines an action by a weak convolution product on the whole of $L^p(\Omega)$, the action commutes with the weak derivative $D_\alpha$ only on the Sobolev spaces which are dense subspaces of $L^p(\Omega)$ closed under the weak derivative. Thus, given a nonlinear map $f \in L^p(\Omega)$, the linear map $\partial f$-the matrix of its partial derivatives-is said to be \emph{regular} at any point where it is invertible. The Inverse Function Theorem, through the invertibility of $\partial f$ the linear part of $f$, posits the invertibility of the nonlinear function $f$ within a neighbourhood of a given (regular) point. This is what constitute the partial dynamical system.
\begin{prop}\normalfont
Linearization preserves invertibility on a functional spaces with the action of the smooth algebra $\mathscr{K}(\Omega)$.
\end{prop}
\begin{proof}
Linearization implies the existence of a linear approximation $f(x) = y + f'(x)x$ for the nonlinear function $f$. The equation $F(x,y) = F_y(x) =  y - f(x)$ from the nonlinear function yields a filter, and its zeros give the points where $f$ is invertible. A net can be derived from the filter given in the form \begin{equation} \mathcal{F}_x := x + \frac{y-f(x)}{Df(a)} \end{equation} which is a form of the Newton approximation formula.
Subsequently, for $y \in \R^n$, we have $\varphi_y(x) := Df(a)^{-1}(y - f(x))$ on $\Omega \subset \R^n$. For every root $x$ of the equation, we have  $f(x) = y \; \implies \; \varphi_y(x) = x$. Hence, what is needed is a space where $\varphi_y(x)$ is a contraction. There is therefore a closed subset $\overline{\Omega} \subset \Omega$ such that the smooth $\mathscr{K}(\Omega)$-action by its net commutes with the partial differential operator $D$.

By the property of the net, the invertibility of $Df_\varepsilon$ follows from that of $Df$; and since $D(f_\varepsilon) = (Df)_\varepsilon$, $f_\varepsilon$ are invertible maps with the limit $f$ also invertible. The commutativity of the smooth $\mathscr{K}(\Omega)$-action with the linear map $Df$ characterizes the Sobolev space $W^{k,p}(\Omega)$. Thus, \cite{HasselblattKatok} sums that if $f \in C^k$ and $g \in C^r$ for some $r < k$. Then $Df(g(y)) \in C^r$ has the inverse $Dg$. Thus $g \in C^{r+1}$.
\end{proof}

Thus, the weak derivatives define the partial dynamical systems of the $L^p$-spaces, the symmetries of which are encoded by smooth algebra $\mathscr{K}(\Omega)$ within the Sobolev Spaces. Alternatively, the Sobolev spaces are spaces of partial dynamical systems of $L^p(\Omega)$. They are given next as space of regular distributions using slice analysis.

\section{Slice Analysis of $\mathscr{K}(\Omega)$-Action}
Noncommutative or nonhomogeneous spaces are generally treated by embedding them in a bigger space with well defined structures, whereby the bad space is seen as the (projective) limit of elements of the bigger space. In \cite{Khalkhali2009}, this approach was  defined  as (equivalence) relations on bigger spaces by forgetting part of their structure. These extra parts lacking in the bad space are the quotient data which constitute the homogeneous or cohomogeneity one space. The Sobolev spaces are treated with slice and cohomogeneity one analysis using this projective approach.

The nets contribute the concept of homogeneity because a net and its limit constitute a unit. It follows, therefore, that the $\mathscr{K}(\Omega)$-action constitutes the homogeneous fibre. Thus, the smooth nets define the embedding $\mathscr{K}(\Omega) : L^p(\Omega) \hookrightarrow \mathcal{D}'(\Omega)$ under the weak product and derivative when they are understood in the sense of  \cite{Marsden68} as a weak limit of 'quantities'. The embedding is used to bring the nets and their limits together as a homogenous unit. The closure and boundedness (or compactness) which follows guarantee the continuity of transformations.

The smoothness property, on the other hand, defines the concept of slices and uniqueness of sections (and vice versa) in the resulting bundle space.
\begin{prop}\normalfont
A Sobolev space $W^{k,p}(\Omega) \subset L^p(\Omega)$ is a $\mathscr{N}(\Omega)$-kernel over the projection $\pi : L^p(\Omega) \to L^p(\Omega)/\mathscr{K}(\Omega)$ and a $\mathscr{N}(\Omega)$-slice in $L^p(\Omega)$; where $\mathscr{N}(\Omega)$ is the set of zero or negligible nets of $C^\infty(\Omega)$ at any point $u \in L^p(\Omega)$.
\end{prop}
\begin{proof}
Applying the slice theorem of \cite{PalaisTerng87}, a Sobolev space $W^{k,p}$ as a closed subspace  of $L^p(\Omega)$ under weak derivative, can be considered a $\mathscr{N}(\Omega)$-kernel. Since $\mathscr{K}(\Omega)$ has the form of locally compact group of transformations of $L^p(\Omega)$ by the proper map $(u_\varepsilon,\varphi_\varepsilon u_\varepsilon) \mapsto (u,\varphi u)$. This implies that $\mathscr{K}(\Omega)$-action is free and proper. So, $L^p(\Omega)$ is a proper $\mathscr{K}(\Omega)$-space because each orbit is closed in $L^p(\Omega)$. Let $\mathscr{N}(\Omega)_u$ be negligible nets at $u \in L^p(\Omega)$, then $\varphi_\varepsilon \mapsto \varphi_\varepsilon u$ is an open map of $\mathscr{K}(\Omega)$ onto $\mathscr{K}(\Omega)(u)$. Hence, given a neighbourhood $W$ of $\mathscr{N}(\Omega)u$ in $\mathscr{K}(\Omega) \subset C^\infty(\Omega)$, it contains an open subset with compact closure in $\mathscr{K}(\Omega)$ acting properly at $u \in L^p(\Omega)$.

Next, given that $S \subset L^p(\Omega)$ is an $\mathscr{N}(\Omega)$-kernel over the $L^p$-projection $\pi : L^p(\Omega) \to L^p(\Omega)/\mathscr{K}(\Omega)$ if there exists an equivariant map $\sigma : \mathscr{K}(\Omega)S \to \mathscr{K}(\Omega)/\mathscr{N}(\Omega)$ such that $\sigma(\mathscr{N}(\Omega)) = S$. With $L^p(\Omega)$ a $\mathscr{K}(\Omega)$-space, and $\mathscr{N}(\Omega)$ a closed subspace/algebra of $\mathscr{K}(\Omega)$, it follows that $S = W^{k,p}(\Omega)$ satisfies the definition while the weak derivative $D$ corresponds to the equivariant map $\sigma$, since $D^{-k}(\mathscr{N}(\Omega)) = W^{k,p}(\Omega)$, and $\mathscr{K}(\Omega)W^{k,p}(\Omega)$ is open in $L^p(\Omega)$.Thus, $W^{k,p}(\Omega)$ is an $\mathscr{N}(\Omega)$-slice in $L^p(\Omega)$.
\end{proof}
\begin{coro}\normalfont
The equivariant map $D : U \to \mathscr{K}(\Omega)$ is a local cross section on the quotient space $\mathscr{K}(\Omega)/\mathscr{N}(\Omega)$.
\end{coro}
\begin{proof}
This follows from the fact that a slice at $u \in L^p(\Omega)$ is a $\mathscr{K}(\Omega)_u$-slice in $L^p(\Omega)$ containing $u$; and an $\mathscr{N}(\Omega)$-slice $W^{k,p}(\Omega)$ defines an equivariant map by the fact that $W^{k,p}(\Omega) = D^{-k}(\mathscr{N}(\Omega))$, where $D^k : \mathscr{K}(\Omega)W^{k,p}(\Omega) \to \mathscr{K}(\Omega)/\mathscr{N}(\Omega)$. This implies that the Sobolev spaces exclude the functions whose weak derivatives vanish at any fixed point $u$ of $L^p(\Omega)$. Since by equivariance $D^k(\varphi_\varepsilon u) = \varphi_\varepsilon D^k u = \varphi_\varepsilon \mathscr{N}(\Omega)_u$, it follows that $D^k$ is determined by $W^{k,p}(\Omega)$. Thus, $D^k \leftrightarrow W^{k,p}(\Omega)$.

Subsequently, on the quotient space $\mathscr{K}(\Omega)/\mathscr{N}(\Omega)$, for a fixed net $\psi_\varepsilon \in \mathscr{K}(\Omega)$, the map $D : W \to \mathscr{K}(\Omega)$, where $W \subset \mathscr{K}(\Omega)/\mathscr{N}(\Omega)$ is an open neighbourhood of $\mathscr{N}(\Omega)$, is a homeomorphism and differentiable, with $D(\mathscr{N}(\Omega)) = e$ and $D(\phi_\varepsilon) = \phi_\varepsilon$ for all $\phi_\varepsilon \in W$. Hence, $D^k$ is a local cross section for $W = W^{k,p}(\Omega)$ a slice at $u \in L^p(\Omega)$.
\end{proof}

Because the weak derivative of a function $f \in L^p(\Omega)$ is unique up to a set of measure zero, $W^{k,p}(\Omega)$ is constituted as spaces of germs of functions or equivalence classes of functions with the same zero sets. This means that the Sobolev spaces $W^{k,p}(\Omega)$ can be given a quotient structure $W^{k,p}(\Omega)/\sim$, where $\sim$ is defined with respect to the convergence of net of measures $\mu_h \to \mu$ or the convergence of nets of measurable functions $f_\varepsilon \to f$; where and $\mu_\varepsilon \to 0$ is a zero net.

The smooth algebra $\mathscr{K}(\Omega) \subset C^\infty(\Omega)$ acts on $W^{k,p}(\Omega)$ by the smooth nets which is a Lie group $\mathscr{K}(\Omega')$ for a compact subset $\Omega' \subset \Omega$, and a Lie pseudogroup action when $\Omega$ is open. Let $H^{k,p}(\Omega)$ denote the closure of the space $C^\infty(\Omega)\cap W^{k,p}(\Omega)$ relative to the $W^{k,p}$-norm. Hence, each equivalence class of $W^{k,p}(\Omega)$ has a fibre action: \begin{equation} \mathscr{K}(\Omega) \times H^{k,p}_o(\Omega) \to H^{k,p}_o(\Omega), \end{equation} shown by the local sections above. So $f \in H^{k,p}_o(\Omega)$ implies existence of a net $(\phi_\varepsilon) \subset C^\infty_o(\Omega)$ such that $\phi_\varepsilon(f) = (f_\varepsilon) \subset H^{k,p}_o(\Omega)$ with $\underset{\varepsilon \to 0}\lim ||f_\varepsilon - f||_{k,p} = 0$.

Since $\varepsilon < dist(x, \partial \Omega)$, it follows that the functions in $H^{k,p}_o(\Omega)$ vanish on the boundary $\partial \Omega$ as $\varepsilon \to 0$. The set of such functions in the $0$-class constitute an ideal of the partial algebra $H^{k,p}_o(\Omega)$ defined by the convolution product of $\mathscr{K}(\Omega)$ on $L^p(\Omega)$ which is compatible with weak derivatives on the Sobolev spaces. With the convergence of every net in $W^{k,p}(\Omega)$, the space is complete with respect to the norm $||\cdot||_{k,p}$, which makes them Banach spaces, with $H^{k,p}(\Omega) = W^{k,p}(\Omega)$. (cf. \cite{Jost2005}, 20.9 and 20.10)

The closure of the spaces $W^{k,p}(\Omega)$ under convolution product and weak derivative as subspaces of $L^p(\Omega)$ and the action of the smooth algebra $\mathscr{K}(\Omega)$ by the nets imply the existence of solutions for every system of partial differential equations defined on a subset $\Omega \subset \R^n$. This follows from the fact that the smooth algebra $\mathscr{K}(\Omega)$ is made up of contractive transformations on a Banach space or manifold. Thus, the convolution product operation on the $L^p$-spaces which defines the action of the smooth algebra $\mathscr{K}(\Omega)$ makes $W^{k,p}(\Omega)$ the spaces of $\mathscr{K}(\Omega)$-generalized quantities.
\begin{thm}\normalfont
The smooth $\mathscr{K}(\Omega)$-action and it invariance under the weak derivative $D^k$ make Sobolev spaces $W^{k,p}(\Omega)$ regular space of generalized $\mathscr{K}(\Omega)$ quantities.
\end{thm}
\begin{proof}
Following Marsden \cite{Marsden68}, the weak $^*$-topology is defined by the nets of functionals $f_\varepsilon \in \mathscr{K}(\Omega)^*$; where $f_\varepsilon \to f \in \mathscr{K}(\Omega)^* \; \iff \; f_\varepsilon(\varphi) \to f(\varphi) \in \R, \; \forall \; \varphi \in \mathscr{K}(\Omega)$, $f_\varepsilon \in W^{k,p}(\Omega)$. This is the result of the weak product and the action of $\mathscr{K}(\Omega)$ on $L^p(\Omega)$. It shows the denseness of $C^\infty$-functions (with compact support) in $L^p(\Omega)$.

The distributional quantities or generalized $\mathscr{K}(\Omega)$ quantities are the $\mathscr{K}(\Omega)$ orbits in $L^p(\Omega)$, which are nets of functionals $f_\varepsilon$ in $\mathscr{K}(\Omega)^* \subset L^p(\Omega) \subset \mathcal{D}'(\Omega)$. Hence, that $f \in \mathscr{K}(\Omega)^*$ is a generalized $\mathscr{K}(\Omega)$ quantity is another way of saying that $f \in \mathscr{K}(\Omega)^*$ defines a weak product (distribution) on $\mathscr{K}(\Omega)$ (or on $C^\infty_c(\R^n)$), and has weak derivative which defines a distribution. The subalgebra (or partial algebra) structure defined on $L^p(\Omega)$ by the action of the smooth algebra $\mathscr{K}(\Omega)$ is due to the fact that the transformations defined by the product and weak derivative are closed in the Sobolev spaces which are embedded in the space of distributions $\mathcal{D}'(\Omega)$.

The Sobolev spaces contain only regular distributions because their definition excludes functions $f$ which cannot be given as local integral of a derivative; that is, whose derivatives $f' = 0$ in the given neighbourhood of a point $x$. Therefore, the action $\mathscr{K}(\Omega) \times W^{k,p}(\Omega) \to W^{k,p}(\Omega)$ is a regularisation and defines transformations of regular distributions. The net approach presents the Sobolev spaces $W^{k,p}(\Omega)$ as regular spaces of generalised quantities of the smooth algebra $\mathscr{K}(\Omega)$. The following definition and proposition complete the proof of the theorem.
\end{proof}
\begin{defx}
The injection $\varphi : \mathfrak{G} \to \mathfrak{G}_c(\pi)^*$ defined on the set $\mathfrak{G}$ of the sections of the fibre bundle $\pi : \mathcal{D}'(\Omega) \to L^p(\Omega)$ (with compact support) are given as \begin{equation} f \mapsto \varphi(f)\cdot g = \int_\Omega \langle f(a), g(a)\rangle d\mu(a). \end{equation} They are called \emph{generalised sections}. Cf. \cite{Marsden68}.
\end{defx}
The generalized sections act on the generalized quantities $\mathscr{K}(\Omega)^*$. This gives the following results.
\begin{prop}\normalfont
The weak derivatives $D^k$ are the generalised sections of a smooth bundle.
\end{prop}
\begin{proof}
By definition, the orbit of the $\mathscr{K}(\Omega)$-action are also nets of integrable functions whose weak derivatives are also of weak derivatives of functions in $W^{k,p}(\Omega)$ (or regular distributions). These constitute the class of regular distributions in $\mathcal{D}'(\Omega)$.

Given any multi-index $\alpha$, and given the differential operator $\partial^\alpha$ and the distribution operator $T$, the weak derivative is denoted $D^\alpha = \partial^\alpha T$, while the distribution defined by a derivative is given as $T_{\partial^\alpha}$, where $T_{\partial^\alpha f} = (-1)^{|\alpha|}(\partial^\alpha T_f)$. define the weak derivative of all orders for a function $\phi \in W^{k,p}(\Omega)$ as $ \langle \partial^\alpha T, \phi \rangle = (-1)^{|\alpha|}\langle T,\partial^\alpha \phi \rangle. $
The definition satisfies $\partial^\alpha \partial^\beta T = \partial^\beta \partial^\alpha T$. Hence, the weak derivative $D^\alpha$ represents a continuous functional on $\mathcal{D}(\Omega)$ since $\phi \in \mathcal{D}(\Omega)$.

The Sobolev spaces are subspaces where the weak derivative of $f \in C^o(\Omega)$ coincides with its classical derivative; that is $\partial^\alpha T_f = T_{\partial^\alpha f}$. Such functions in $L^p(\Omega)$ are said to define regular distributions. The differential operator $\partial^\alpha$ is said to be both $\mathscr{K}(\Omega)$ and distribution $T$-invariant on such subspaces of distributions. The Sobolev spaces $W^{k,p}(\Omega) \subset \mathcal{D}'(\Omega)$ are therefore subspaces of regular distributions on which the classical and weak differential operator $\partial^\alpha$ are the same and commutes with the distributional operator $T$; that is $\partial^\alpha T_f = T_{\partial^\alpha f}$. The weak derivatives $D^\alpha$, for $\alpha \leq k$, are sections of the smooth bundle $\pi : \mathcal{D}'(\Omega) \to L^p(\Omega)$. The smooth algebra action codes this dynamical symmetry on the space of distributions.
\end{proof}
\begin{prop}\normalfont
The action of the smoothing algebra $\mathscr{K}(\Omega)$ by its net attaches a smooth fibre to the regular generalized functions (or distributions) in $L^p(\Omega)$.
\end{prop}
\begin{proof}
Considered as closed and embedded subspaces of distributions $D'(\Omega)$, the Sobolev spaces $W^{k,p}(\Omega)$ are subalgebras of the algebra of distributions. In the generalized setting, as in \cite{KunzingerKonjik}, given a distributional vector field $X \in D'(M,TM)$ on a differentiable manifold $M$, which is the limit of a net of smooth vector fields $(X_\varepsilon)_\varepsilon$ possessing complete flows $\varphi_\varepsilon(t,.)$; $X$ is shown to possess a measurable flow $\phi_t = \varphi(t,.)$ whenever (i) $X_\varepsilon \to X \; \implies \; \varphi_\varepsilon(t,.) \to \varphi(t,.)$, almost everywhere on $M$ for all $t$ and $\phi_t$ is measurable; (ii) For each $t \in \R$ and each $C \subset \subset M$ there exists $\varepsilon_o \in I$ and $K \subset \subset M$ with $C \subseteq K$ such that $\varphi_\varepsilon(t,C) \subseteq K$ for all $\varepsilon$. That is, $\varphi_\varepsilon(t,.)$ is compactly supported in $M$ or $c$-bounded.
\end{proof}
\begin{coro}\normalfont
The convergence of smooth functions in $W^{k,p}(\Omega)$ implies the possibility of projection on vector fields.
\end{coro}
\begin{proof}
The idea of associating generalized vector fields to distributional vector fields is similar to that of project-ability of vector fields. The projection of a generalized vector field $u \in \mathcal{G}(M)$ is its distributional shadow $\omega \in D'(M)$, satisfying the equality $\displaystyle \underset{\varepsilon \to 0}\lim \int_M u_\varepsilon v = \langle w,v \rangle$, where $v$ is a compactly supported $1$-density on $M$. Thus, the relation $\sim$ of association generalizes distributional equality or weak equality to the generalized space $\mathcal{G}(M)$. But on the Sobelev spaces, it means convergence of nets of smooth functions or distributions.
\end{proof}
\begin{coro}\normalfont
The solution to the initial value problem $\dot{x} = kx(t), \; x(0) = x_o$, is the limit of the generalized solution of a more generalized problem $\dot{x}(t) = X(x(t)), \; x(t_o) = x_o$, where $X \in D'(M,TM)$ is a distributional vector field as defined in \cite{KunzingerKonjik}.
\end{coro}
\begin{proof}
In generalized algebra, a moderate net $Y = (Y_\varepsilon)_{\varepsilon \in I} \in \mathcal{G}^1_o(M)$ is associated to a smooth distributional vector field $X \in D'(M, TM)$. Thus, the map $\mathcal{G}^1_o(M) \to D'(M,TM); \; Y \to X$ is a smooth projection. Since the generalized element $Y$ is a smooth net of moderate fields, which generate the flows $\varphi_\varepsilon(t,.)$, these flows are also associated with the flow $\varphi(t,.)$ generated by $X$. This then means that $X$ is a distribution if there is a moderate generalized field $Y = (Y_\varepsilon)_\varepsilon$ associated with $X$. This is what we have established as a section $D^k$ in the case of the Sobolev spaces and their geometry.
\end{proof}

\section{Conclusion}
Because Sobolev spaces $W^{k,p}(\Omega)$ are dense in $L^p(\Omega)$, they are treated as embedded subspaces of $L^p(\Omega)$. Hence, the convergence of its elements $f_\varepsilon \to f$ is used to study the $L^p$-spaces via their distributions. Also, given that the Sobolev spaces are closed under weak derivatives, they form closed embedded submanifolds of distributions $D'(\Omega)$, on which weakly differentiable maps are contractions.

Theorem 6.2 of \cite{Marsden68} modified by \cite{Kunzinger2003} connects this to the definition of partial dynamical systems, which is formally defined by Exel \cite{Exel2017}. Thus, given $X \in D'(M,TM)$ a vector field with measurable flow $\varphi(t,.)$, the flow property is said to hold in the sense that $\varphi_{t+s} = \varphi_t \circ \varphi_s$ almost everywhere on $M$, for all $s,t \in \R$. From this property, it follows that the smooth algebra action establishes the structure of distribution and the regularization of flows in the weakly differentiable manifolds of $L^p(\Omega)$. The Sobolev spaces are established in this work to be the stable spaces of the partial dynamical systems defined by nets of the smooth algebra $\mathscr{K}(\Omega)$ on the $L^p$-spaces.\\

\end{document}